# MEAN SQUARED ERROR OF EMPIRICAL PREDICTOR

By Kalyan Das[1], Jiming Jiang[2] and J. N. K. Rao[1]

*Calcutta University, University of California, Davis and Carleton University*

The term "empirical predictor" refers to a two-stage predictor of a linear combination of fixed and random effects. In the first stage, a predictor is obtained but it involves unknown parameters; thus, in the second stage, the unknown parameters are replaced by their estimators. In this paper, we consider mean squared errors (MSE) of empirical predictors under a general setup, where ML or REML estimators are used for the second stage. We obtain second-order approximation to the MSE as well as an estimator of the MSE correct to the same order. The general results are applied to mixed linear models to obtain a second-order approximation to the MSE of the empirical best linear unbiased predictor (EBLUP) of a linear mixed effect and an estimator of the MSE of EBLUP whose bias is correct to second order. The general mixed linear model includes the mixed ANOVA model and the longitudinal model as special cases.

**1. Introduction.** We consider a general linear mixed model of the form

$$(1.1) \qquad y = X\beta + Zv + e,$$

where $y$ is an $n \times 1$ vector of sample observations, $X$ and $Z$ are known matrices, $\beta$ is a $p \times 1$ vector of unknown parameters (fixed effects) and $v$ and $e$ are distributed independently with means 0 and covariance matrices $G$ and $R$, respectively, depending on some unknown vector of parameters $\sigma$. We assume that $p$ is fixed and $X$ is of full rank $p$ $(< n)$. Note that $\mathrm{cov}(y) = \Sigma = R + ZGZ'$.

Problems involving multiple sources of random variation are often modeled as special cases of (1.1). For example, in the well-known ANOVA

Received August 2001; revised March 2003.

[1]Supported in part by a National Sciences and Engineering Research Council of Canada grant.

[2]Supported in part by NSF Grant SES-99-78101.

model we partition $Z$ as $Z = (Z_1, \ldots, Z_q)$ and $v = (v_1', \ldots, v_q')'$, where $Z_i$ is $n \times r_i$, $v_i$ is $r_i \times 1$, $i = 1, \ldots, q$, and $v_1, \ldots, v_q$ are mutually independent with $v_i \sim N(0, \sigma_i I_{r_i})$ and $e \sim N(0, \sigma_0 I_n)$. (For notational convenience we use $\sigma_i$ rather than the customary $\sigma_i^2$ to denote the $i$th variance component.) Note that the ANOVA model is a special case of (1.1) with $R = \sigma_0 I_n$, $G = \mathrm{diag}(\sigma_1 I_{r_1}, \ldots, \sigma_q I_{r_q})$ and $\sigma = (\sigma_0, \sigma_1, \ldots, \sigma_q)'$. The (dispersion) parameter space under the ANOVA model is $\Theta = \{\sigma : \sigma_i \geq 0, \ i = 0, 1, \ldots, q\}$. The well-known "longitudinal" random effects model [Laird and Ware (1982)] is also a special case of (1.1). In this case $y = (y_1', \ldots, y_t')'$ with

$$(1.2) \qquad y_i = X_i \beta + Z_i v_i + e_i, \qquad i = 1, \ldots, t,$$

where $y_i$ is $n_i \times 1$, $X_i$ is $n_i \times p$ and $Z_i$ is $n_i \times r_i$. It is assumed that the $y_i$'s are independent, $\mathrm{cov}(v_i) = G_i$, $\mathrm{cov}(e_i) = R_i$, where $G_i$ and $R_i$ depend on $\sigma$, and $v_i$ and $e_i$ are independent. It follows that $\Sigma = \mathrm{diag}(\Sigma_1, \ldots, \Sigma_t)$ with $\Sigma_i = \mathrm{cov}(y_i) = \Sigma_i = R_i + Z_i G_i Z_i'$. [Note that the longitudinal model (1.2) is not a special case of the ANOVA model and vice versa.] The well-known Fay–Herriot (1979) model widely used in small area estimation is a special case of the longitudinal model. The (dispersion) parameter space under the longitudinal model is $\Theta = \{\sigma : \Sigma_i$ is nonnegative definite, $i = 1, \ldots, t\}$.

Estimation of linear combinations of $\beta$ and realized $v$ from (1.1), say $\mu = l'\beta + m'v$, for specified vectors of constants $l$ and $m$ is of considerable interest in many practical applications, for example, the estimation of quality index, longitudinal studies, the selection index in quantitative genetics, plant varietal trials and small area estimation [Robinson (1991)]. Henderson (1975) obtained the best linear unbiased predictor (BLUP) of $\mu$ under model (1.1) as

$$
\begin{aligned}
(1.3) \quad t(\sigma) &= t(\sigma, y) \\
&= l'\tilde{\beta} + m'\tilde{v} = l'\tilde{\beta} + s(\sigma)'(y - X\tilde{\beta}),
\end{aligned}
$$

where

$$\tilde{\beta} = \tilde{\beta}(\sigma) = (X'\Sigma^{-1}X)^{-1}X'\Sigma^{-1}y$$

is the generalized least squares estimator, or best linear unbiased estimator (BLUE), of $\beta$, $\tilde{v} = \tilde{v}(\sigma) = GZ'\Sigma^{-1}(y - X\tilde{\beta})$ and $s(\sigma) = \Sigma^{-1}ZGm$.

The BLUP estimator (1.3) is unbiased in the sense of $\mathrm{E}[t(\sigma, y) - \mu] = 0$ under (1.1). The mean squared error (MSE) of $t(\sigma)$ is given by

$$(1.4) \qquad \mathrm{MSE}[t(\sigma)] = \mathrm{E}[t(\sigma) - \mu]^2 = g_1(\sigma) + g_2(\sigma),$$

where

$$g_1(\sigma) = m'(G - GZ'\Sigma^{-1}ZG)m$$

and

$$g_2(\sigma) = [l - X's(\sigma)]'(X'\Sigma^{-1}X)^{-1}[l - X's(\sigma)];$$

see Henderson (1975). Results (1.3) and (1.4) do not require normality of random effects $v$ and $e$.

The BLUP estimator $t(\sigma)$ depends on the dispersion parameters $\sigma$, which are unknown in practice. It is therefore necessary to replace $\sigma$ by a consistent estimator $\hat{\sigma}$ to obtain a two-stage estimator or empirical BLUP (EBLUP) given by $t(\hat{\sigma})$. Methods of estimating $\sigma$ include maximum likelihood (ML) and restricted maximum likelihood (REML) under normality, the method of fitting-of-constants and minimum norm quadratic unbiased estimation (MINQUE) without the normality assumption; see Searle, Casella and Mc-Culloch (1992). The resulting estimators $\hat{\sigma}$ are even and translation invariant, that is, $\hat{\sigma}(y) = \hat{\sigma}(-y)$ for all $y$ and $\hat{\sigma}(y + X\beta) = \hat{\sigma}(y)$ for all $y$ and $\beta$. Jiang (1996) proved that ML and REML estimators $\hat{\sigma}$ obtained under normality remain consistent without the normality assumption.

Kackar and Harville (1981) showed that the EBLUP $t(\hat{\sigma})$ remains unbiased if $\hat{\sigma}$ is even and translation invariant. This result holds provided that $\mathrm{E}[t(\hat{\sigma})]$ is finite and $v$ and $e$ are symmetrically distributed (not necessarily normal). In particular, the two-stage estimator $\hat{\beta} = \tilde{\beta}(\hat{\sigma})$ is unbiased for $\beta$.

Kenward and Roger (1997) studied inference for the fixed effects, $\beta$, under a general Gaussian linear mixed model $y \sim N(X\beta, \Sigma)$ with a structured covariance matrix $\Sigma = \Sigma(\sigma)$ depending on some parameter $\sigma$. They used the REML estimator of $\beta$, namely the two-stage estimator $\hat{\beta} = \tilde{\beta}(\hat{\sigma})$, where $\hat{\sigma}$ is the REML estimator of $\sigma$. A naive estimator of $\mathrm{cov}(\hat{\beta})$ that ignores the variability in $\hat{\sigma}$ is given by $[X'\Sigma^{-1}(\hat{\sigma})X]^{-1}$. Kenward and Roger (1997) derived a bias-adjusted estimator of $\mathrm{cov}(\hat{\beta})$ and used it to derive a scaled Wald statistic, together with an F approximation to its distribution. The F approximation performed well in simulations under a range of small sample settings. Kenward and Roger (1997) did not study the precise order of the bias of the adjusted estimator.

Kackar and Harville (1981) studied approximation to the MSE of EBLUP $t(\hat{\sigma})$, assuming normality of the random effects $v$ and errors $e$ in the model (1.1). They showed that

$$(1.5) \qquad \mathrm{MSE}[t(\hat{\sigma})] = \mathrm{MSE}[t(\sigma)] + \mathrm{E}[t(\hat{\sigma}) - t(\sigma)]^2$$

for any even and translation invariant estimator $\hat{\sigma}$, provided that $\mathrm{MSE}[t(\hat{\sigma})]$ is finite. It should be pointed out that, under very mild conditions, $\mathrm{E}[t(\hat{\sigma})]$ and $\mathrm{MSE}[t(\hat{\sigma})]$ are, in fact, finite [see Jiang (2000)]. It is customary among practitioners to ignore the variability associated with $\hat{\sigma}$ and use the following naive estimator of $\mathrm{MSE}[t(\hat{\sigma})]$:

$$(1.6) \qquad \mathrm{mse}_{\mathrm{N}}[t(\hat{\sigma})] = g_1(\hat{\sigma}) + g_2(\hat{\sigma}).$$

However, it follows from (1.4) and (1.5) that (1.6) can lead to significant underestimation. Therefore, it is practically important to obtain approximately unbiased estimators of MSE that reflect the true variability associated with the EBLUP estimators. This becomes particularly important when large sums of funds are involved. For example, over \$7 billion dollars of funds are allocated annually on the basis of EBLUP estimators of school-age children in poverty at the county and school district levels [National Research Council (2000)].

Kackar and Harville (1984) gave an approximation to $\mathrm{MSE}[t(\hat{\sigma})]$ under the general model (1.1), taking account of the variability in $\hat{\sigma}$, and proposed an estimator of $\mathrm{MSE}[t(\hat{\sigma})]$ based on this approximation. However, the approximation is somewhat heuristic, and the accuracy of the approximation and the associated estimator of $\mathrm{MSE}[t(\hat{\sigma})]$ was not studied. Prasad and Rao (1990) studied the accuracy of a second-order approximation to $\mathrm{MSE}[t(\hat{\sigma})]$ for two important special cases of the longitudinal model (1.2): (i) the well-known Fay–Herriot model (2.15) studied in Section 2.3 and (ii) the nested error linear regression model given by (1.2) with $Z_i = \mathbf{1}_{n_i}$, a scalar $v_i$ with $\mathrm{var}(v_i) = \sigma_1$ and $\mathrm{cov}(e_i) = \sigma_0 I_{n_i}$, $i = 1, \ldots, t$. In the context of small area estimation $n_i$ is the sample size in the $i$th area and $t$ is the number of small areas. The nested error model may also be regarded as a special case of the ANOVA model with a single source of variation ($q = 1$), $G = \sigma_1 I_t$ and $R = \sigma_0 I_n$. Using the method of fitting-of-constants estimator $\hat{\sigma}$, Prasad and Rao (1990) showed that, for large $t$,

$$(1.7) \qquad \mathrm{E}[t(\hat{\sigma}) - t(\sigma)]^2 = g_3(\sigma) + o(t^{-1}),$$

where $g_3(\sigma)$ depends on $\mathrm{cov}(\hat{\sigma})$. This leads to a second-order approximation

$$(1.8) \qquad \mathrm{MSE}_{\mathrm{a}}[t(\hat{\sigma})] = g_1(\sigma) + g_2(\sigma) + g_3(\sigma).$$

The approximation is accurate to terms $o(t^{-1})$, that is, the neglected terms are $o(t^{-1})$. The $g_3(\sigma)$ term is computationally simpler compared to an asymptotically equivalent term obtained from Kackar and Harville's approximation. Prasad and Rao (1990) also obtained an estimator of $\mathrm{MSE}[t(\hat{\sigma})]$ given by

$$(1.9) \qquad \mathrm{mse}_{\mathrm{PR}}[t(\hat{\sigma})] = g_1(\hat{\sigma}) + g_2(\hat{\sigma}) + 2g_3(\hat{\sigma}).$$

The estimator (1.9) is approximately unbiased in the sense that its bias is $o(t^{-1})$. Kackar and Harville (1984) proposed an alternative estimator given by

$$(1.10) \qquad \mathrm{mse}_{\mathrm{KH}}[t(\hat{\sigma})] = g_1(\hat{\sigma}) + g_2(\hat{\sigma}) + g_3^*(\hat{\sigma})$$

for any even and translation invariant estimator $\hat{\sigma}$. The bias of (1.10) is $O(t^{-1})$; that is, it is not approximately unbiased to terms $o(t^{-1})$.

Harville and Jeske (1992) studied various MSE estimators under the ANOVA model with a single source of random variation ($q = 1$) and the REML estimator of $\sigma$, including an estimator of the form (1.9). They referred to the latter estimator as the Prasad–Rao (P–R) estimator. They appealed to Prasad–Rao's asymptotic results for its justification, but the latter results have been justified only for the special cases (i) and (ii). They also conducted a limited simulation study based on the simple one-way random effects model, $y_{ij} = \beta + v_i + e_{ij}$, $j = 1, \ldots, n_i$, $i = 1, \ldots, t$, using a small balanced design ($t = 9$, $n_i = 2$ for all $i$), a small unbalanced design ($t = 9$, $n_1 = \cdots = n_8 = 1$, $n_9 = 10$) and a large unbalanced design ($t = 21$, $n_1 = \cdots = n_{20} = 1$, $n_{21} = 50$). The objective was to estimate the mean $\mu = \beta + v_1$. Simulation results indicated that the P–R estimator performs well when $\gamma = \sigma_1/\sigma_0$ is not small, but it can lead to substantial overestimation for small values of $\gamma$ closer to the lower bound of 0, especially for small $t$ ($= 9$). Two partially Bayes estimators perform better than the P–R estimator when $\gamma$ is close to 0.

Datta and Lahiri (2000) extended Prasad and Rao's (1990) results to the general longitudinal model (1.2) with covariance matrices $R_i$ and $G_i$ having linear structures of the form $R_i = \sum_{j=0}^{q} \sigma_j R_{ij} R'_{ij}$ and $G_i = \sum_{j=0}^{q} \sigma_j G_{ij} G'_{ij}$, where $\sigma_0 = 1$, $R_{ij}$ and $G_{ij}$ ($i = 1, \ldots, t$; $j = 0, 1, \ldots, q$) are known matrices with uniformly bounded elements such that $R_i$ and $G_i$ are positive definite matrices for $i = 1, \ldots, t$. They studied ML and REML estimators of $\sigma = (\sigma_1, \ldots, \sigma_q)'$ and showed that an estimator of MSE[$t(\hat\sigma)$] of the form (1.9) is approximately unbiased when the REML estimator of $\sigma$ is used but not when the ML estimator is used. In the latter case an additional term involving the bias of the ML estimator $\hat\sigma$ is needed for getting an approximately unbiased MSE estimator. Datta and Lahiri (2000) also gave explicit formulas under ML and REML for the two special cases (i) and (ii) studied by Prasad and Rao (1990). The underlying proof of Datta and Lahiri (2000), however, is not rigorous.

The main purpose of our paper is to study the general linear mixed model (1.1) that covers the ANOVA model as well as the longitudinal model and derive a second-order approximation to MSE of EBLUP $t(\hat\sigma)$ under REML and ML estimation of the variance components parameters $\sigma$. We also derive approximately unbiased estimators of MSE[$t(\hat\sigma)$] and specify the precise order of the neglected terms. For ANOVA models with multiple sources of random variation, the components of $\hat\sigma$ may have different convergence rates [Miller (1977) and Jiang (1996)]. As a result, rigorous proofs are quite technical and long. We have therefore only sketched the technical details in Section 5 of our paper, but the detailed proofs are available at the web site address given in Section 5.

The remaining sections of the paper are organized as follows. In Section 2, we first present a general asymptotic representation of $\hat\sigma - \sigma$, where $\hat{\sigma}$ is obtained as a solution of "score" equations of the form $\partial l(\sigma)/\partial\sigma = 0$, and $\sigma$ represents the true value of parameter vector. Normality assumption is not needed for this asymptotic representation. We then verify that the conditions underlying this representation are satisfied by solutions to the ML and REML score equations belonging to a parameter space $\Theta$ under the ANOVA model and normality. As another example, we show that the conditions are satisfied by the ML and REML estimators under the Fay–Herriot model and normality. In Section 3 we obtain a second-order approximation to $\mathrm{MSE}[t(\hat{\sigma})]$ under normality. The second-order approximation is then spelled out under the ANOVA model using ML and REML estimators $\hat{\sigma}$. We also verify the underlying key conditions for the special cases of the balanced ANOVA model and two special cases of the longitudinal model: the Fay–Herriot model and the nested error regression model. Section 4 gives an estimator of $\mathrm{MSE}[t(\hat{\sigma})]$ correct to second order. The MSE estimator is then spelled out under the ANOVA model and the longitudinal model using ML and REML estimators $\hat{\sigma}$. Technical details are sketched in Section 5.

**2. Asymptotic representation of $\hat{\boldsymbol{\sigma}} - \boldsymbol{\sigma}$.** Throughout the rest of this paper, $\sigma$ represents the true parameter vector in places where there is no confusion; expressions such as $\partial l(\hat{\sigma})/\partial\sigma$ mean derivative with respect to $\sigma$ evaluated at $\tilde{\sigma}$; in expressions such as $\mathrm{E}[\partial^2 l(\sigma)/\partial\sigma^2]$, the expectation is taken at the true $\sigma$, and the function inside $\mathrm{E}(\cdot)$ is also evaluated at the true $\sigma$. We first obtain an asymptotic representation of $\hat{\sigma} - \sigma$, where $\hat{\sigma}$ is obtained as a solution to "score" equations of the form

$$(2.1) \qquad \frac{\partial l(\sigma)}{\partial\sigma} = 0$$

and then apply the general theory to the ANOVA model with REML and ML estimation of $\sigma$. In Section 5.1 we sketch the proof of the asymptotic representation. Here $l(\sigma)$ may correspond to the logarithm of the restricted likelihood $l_{\mathrm{R}}(\sigma)$, or the profile loglikelihood $l_{\mathrm{P}}(\sigma)$ under model (1.1) and normality of $v$ and $e$.

THEOREM 2.1. *Suppose that:*

(i) *$l(\sigma) = l(\sigma, y)$ is three times continuously differentiable with respect to $\sigma = (\sigma_1, \ldots, \sigma_s)'$, where $y = (y_1, \ldots, y_n)'$;*

(ii) *the true $\sigma \in \Theta^{\mathrm{o}}$, the interior of $\Theta$;*

(iii)

$$(2.2) \qquad -\infty < \limsup_{n\to\infty} \lambda_{\max}(D^{-1}AD^{-1}) < 0,$$

*where $\lambda_{\max}$ means the largest eigenvalue, $A = \mathrm{E}[\partial^2 l(\sigma)/\partial\sigma^2]$ and $D = \mathrm{diag}(d_1, \ldots, d_s)$ with $d_i > 0$, $1 \le i \le s$, such that $d_* = \min_{1 \le i \le s} d_i \to \infty$ as $n \to \infty$; and*

(iv) *the gth moments of the following are bounded (g > 0):*

$$\frac{1}{d_i}\left|\frac{\partial l(\sigma)}{\partial \sigma_i}\right|, \qquad \frac{1}{\sqrt{d_i d_j}}\left|\frac{\partial^2 l(\sigma)}{\partial \sigma_i \partial \sigma_j} - \mathrm{E}\left[\frac{\partial^2 l(\sigma)}{\partial \sigma_i \partial \sigma_j}\right]\right|, \qquad \frac{d_*}{d_i d_j d_k}M_{ijk}(y),$$

$$1 \leq i,j,k \leq s,$$

*where* $M_{ijk}(y) = \sup_{\tilde{\sigma} \in S_\delta(\sigma)} |\partial^3 l(\tilde{\sigma})/\partial \sigma_i \partial \sigma_j \partial \sigma_k|$ *with* $S_\delta(\sigma) = \{\tilde{\sigma} : |\tilde{\sigma}_i - \sigma_i| \leq \delta d_*/d_i, 1 \leq i \leq s\}$ *for some* $\delta > 0$. *Then* (1) *a* $\hat{\sigma}$ *exists such that for any* $0 < \rho < 1$ *there is a set* $\mathcal{B}$ *satisfying for large* $n$ *and on* $\mathcal{B}$, $\hat{\sigma} \in \Theta$, $\partial l(\hat{\sigma})/\partial \sigma = 0$, $|D(\hat{\sigma} - \sigma)| < d_*^{1-\rho}$ *and*

$$\hat{\sigma} = \sigma - A^{-1}a + r, \tag{2.3}$$

*where* $a = \partial l(\sigma)/\partial \sigma$ *and* $|r| \leq d_*^{-2\rho}u_*$ *with* $\mathrm{E}(u_*^g)$ *bounded, and* (2) $\mathrm{P}(\mathcal{B}^c) \leq cd_*^{-\tau g}$, *where* $\tau = (1/4) \wedge (1-\rho)$ *and* $c$ *is a constant.*

Note that Theorem 2.1 states that the solution to (2.1), $\hat{\sigma}$, exists and lies in the parameter space with probability tending to 1 and gives the convergence rate of $\hat{\sigma}$ to $\sigma$ as well as the asymptotic representation (2.3), assuming that the true vector $\sigma$ belongs to the interior of the parameter space $\Theta$.

2.1. *REML estimation under the ANOVA model.* The ANOVA model is given by

$$y = X\beta + \sum_{i=1}^{q} Z_i v_i + e. \tag{2.4}$$

The restricted loglikelihood under the ANOVA model with normality of $v$ and $e$ has the form

$$l_{\mathrm{R}}(\sigma) = c - (1/2)[\log(|T'\Sigma T|) + y'Py], \tag{2.5}$$

where $\sigma = (\sigma_0, \sigma_1, \ldots, \sigma_q)'$, $c$ is a constant, $|T'\Sigma T|$ is the determinant of $T'\Sigma T$, $T$ is any $n \times (n-p)$ matrix such that $\mathrm{rank}(T) = n - p$ and $T'X = 0$,

$$
\begin{aligned}
P &= T(T'\Sigma T)^{-1}T' \\
&= \Sigma^{-1} - \Sigma^{-1}X(X'\Sigma^{-1}X)^{-1}X'\Sigma^{-1},
\end{aligned}
\tag{2.6}
$$

and $\Sigma = \sum_{i=0}^{q} \sigma_i V_i$ with $V_0 = I_n$ and $V_i = Z_i Z_i'$, $i \geq 1$ [e.g., Searle, Casella and McCulloch (1992), page 451]. The REML estimator of $\sigma$ is a solution to (2.1) with $l(\sigma) = l_{\mathrm{R}}(\sigma)$. Section 5.2 sketches the proof that shows the conditions of Theorem 2.1 are satisfied, provided that the same conditions under which REML estimators are consistent are satisfied [Jiang (1996)]. The actual proof is somewhat lengthy and uses results on moments of quadratic forms in normal variables and $d_i = \|Z_i'PZ_i\|_2$ for the ANOVA model, where $\|B\|_2 =$ $[\mathrm{tr}(B'B)]^{1/2}$ for a matrix $B$. A quadratic form in $y$ or $u = y - X\beta$ appears in the formulas for the first derivatives of $l_{\mathrm{R}}$,

$$(2.7) \quad \begin{aligned} \frac{\partial l_{\mathrm{R}}(\sigma)}{\partial \sigma_i} &= \frac{1}{2}[y'PV_iPy - \mathrm{tr}(PV_i)] \\ &= \frac{1}{2}[u'PV_iPu - \mathrm{tr}(PV_i)], \qquad 0 \le i \le q. \end{aligned}$$

Note that $u \sim N(0, \Sigma)$. Similarly, the second and third derivatives involve quadratic forms in $u$; see Section 5.2.

2.2. *ML estimation under the ANOVA model.* The (unrestricted) log-likelihood has the form

$$(2.8) \qquad l(\beta, \sigma) = c - \tfrac{1}{2}[\log(|\Sigma|) + (y - X\beta)'\Sigma^{-1}(y - X\beta)],$$

where $c$ is a constant. We have

$$(2.9) \qquad \frac{\partial l(\beta, \sigma)}{\partial \beta} = X'\Sigma^{-1}y - X'\Sigma^{-1}X\beta,$$

$$(2.10) \qquad \frac{\partial l(\beta, \sigma)}{\partial \sigma_i} = \frac{1}{2}[(y - X\beta)'\Sigma^{-1}V_i\Sigma^{-1}(y - X\beta) - \mathrm{tr}(\Sigma^{-1}V_i)],$$
$$0 \le i \le q.$$

Solving $\partial l(\beta, \sigma)/\partial \beta = 0$ for $\beta$, we obtain from (2.9) $\tilde{\beta}(\sigma) = (X'\Sigma^{-1}X)^{-1}X' \times \Sigma^{-1}y$. Substituting $\tilde{\beta}(\sigma)$ for $\beta$ in (2.8), and using (2.6), we obtain the profile loglikelihood

$$(2.11) \qquad l_{\mathrm{P}}(\sigma) = c - \tfrac{1}{2}[\log(|\Sigma|) + y'Py].$$

It now follows that the MLE of $\sigma$ is the solution of the equation

$$(2.12) \qquad \frac{\partial l_{\mathrm{P}}(\sigma)}{\partial \sigma} = 0.$$

Note that $l_{\mathrm{P}}(\sigma)$ is not a loglikelihood function, but Theorem 2.1 does not require $l(\sigma)$ to be a loglikelihood function, so we can take $l(\sigma) = l_{\mathrm{P}}(\sigma)$ in Theorem 2.1. Section 5.3 sketches the proof that shows the conditions of Theorem 2.1 are satisfied with the same $d_i$ as in the REML case and the same set of conditions, provided $p$, the dimension of $\beta$, is bounded as $n$ increases. Again, quadratic forms appear in the formulas for the first derivatives of $l_{\mathrm{P}}(\sigma)$:

$$(2.13) \quad \begin{aligned} \frac{\partial l_{\mathrm{P}}(\sigma)}{\partial \sigma_i} &= \frac{1}{2}[y'PV_iPy - \mathrm{tr}(\Sigma^{-1}V_i)] \\ &= \frac{1}{2}[u'PV_iPu - \mathrm{tr}(\Sigma^{-1}V_i)], \qquad 0 \le i \le q. \end{aligned}$$

Similarly, the second and third derivatives involve quadratic forms in $u$; see Section 5.3.

2.3. *The Fay–Herriot model.* In Sections 2.1 and 2.2 we considered ML and REML estimations under the ANOVA model. Now we consider a different case, the Fay–Herriot model [Fay and Herriot (1979)]. This model has been considered by many authors; see Ghosh and Rao (1994) for a review and extensions.

Suppose that $y_i$ is a scalar random variable such that

$$(2.14) \qquad y_i = x_i'\beta + v_i + e_i, \qquad i = 1, \ldots, t,$$

where the $v_i$'s are i.i.d. $\sim N(0, \sigma)$, $e_i$'s are independent such that $e_i \sim N(0, \phi_i)$ with known $\phi_i$, and $v_i$'s are independent of $e_i$'s. Furthermore, $x_i$ is a known $p \times 1$ vector of covariates, and $\beta$ is an unknown vector of regression coefficients. In the context of small area estimation, $y_i$ denotes a survey estimate of the $i$th area mean $\mu_i$ and $e_i$ denotes the sampling error with known sampling variance, $\mathrm{var}(e_i) = \phi_i$. Furthermore, $\mu_i$ is modelled as $\mu_i = x_i'\beta + v_i$ with model errors $v_i$.

Note that model (2.14) is not a special case of the ANOVA model. In fact, it may be considered as a special case of the longitudinal model introduced in Section 1. Model (2.14) may be written in matrix form as

$$(2.15) \qquad y = X\beta + v + e,$$

where $y = (y_1, \ldots, y_t)'$, $v = (v_1, \ldots, v_t)' \sim N(0, \sigma I_t)$, $e = (e_1, \ldots, e_t)' \sim N(0, \Phi)$ with $\Phi = \mathrm{diag}(\phi_1, \ldots, \phi_t)$, $X$ is $t \times p$ with $i$th row $x_i'$, and $v$, $e$ are independent.

Now consider REML and ML estimations under the Fay–Herriot model (2.15). In Section 5.6 we sketch the proofs that show the conditions of Theorem 2.1 are satisfied if REML or ML estimators of $\sigma$ are used, provided that $\sigma$ is positive and the $\phi_i$'s are bounded.

**3. MSE approximation.** We now obtain a second-order approximation to the MSE of EBLUP $t(\hat{\sigma})$. Under normality the MSE of $t(\hat{\sigma})$ satisfies (1.5), that is,

$$(3.1) \qquad \mathrm{MSE}[t(\hat{\sigma})] = \mathrm{MSE}[t(\sigma)] + \mathrm{E}[t(\hat{\sigma}) - t(\sigma)]^2,$$

where $\mathrm{MSE}[t(\sigma)]$ is obtained from (1.4). It remains to approximate the last term on the right-hand side of (3.1).

THEOREM 3.1. *Suppose that the conditions of Theorem* 2.1 *are satisfied. Furthermore, suppose that* $t(\sigma)$ *can be expressed as*

$$(3.2) \qquad t(\sigma) = \sum_{k=1}^{K} \lambda_k(\sigma) W_k(y),$$

*where* $K = O(d_*^\alpha)$ *for some* $\alpha \geq 0$, *and the following terms are bounded for some* $b > 2$ *and* $\delta > 0$:

(i)
$$\max_{1 \le k \le K} \mathrm{E}|W_k(y)|^b,$$

(ii)
$$\max_{1 \le k \le K} \sup_{\sigma} |\lambda_k(\sigma)|,$$

(iii)
$$\sum_{k=1}^{K} \left| \frac{\partial \lambda_k(\sigma)}{\partial \sigma} \right|,$$

(iv)
$$\sum_{k=1}^{K} \sup_{|\tilde{\sigma} - \sigma| \le \delta} \left\| \frac{\partial^2 \lambda_k(\tilde{\sigma})}{\partial \sigma^2} \right\|,$$

where $\|B\| = [\lambda_{\max}(B'B)]^{1/2}$ is the spectral norm of a matrix $B$. If $g > 8(1 + \alpha)(1 - 2/b)^{-1}$, then

(3.3)            $$\mathrm{E}[t(\hat{\sigma}) - t(\sigma)]^2 = \mathrm{E}(h'A^{-1}a)^2 + o(d_*^{-2}),$$

where $h = \partial t(\sigma)/\partial \sigma$, $A = \mathrm{E}[\partial^2 l(\sigma)/\partial \sigma]$, and $a = \partial l(\sigma)\, \partial \sigma$.

A sketch of the proof of Theorem 3.1 is given in Section 5.4. Note that normality is not required in Theorem 3.1. On the other hand, Theorem 3.1 requires that the predictor $t(\sigma)$ have the form (3.2). In the next two subsections we show that this condition holds for balanced ANOVA models as well as for two longitudinal models. It is possible to replace (3.2) by some moment conditions on $t(\sigma)$ and its first and second derivatives, provided that one considers instead a truncated version of $\hat{\sigma}$, which is defined as $\hat{\sigma}$ if $|\hat{\sigma}| \le L_n$, and is $\sigma^*$ otherwise with $\sigma^*$ being a known vector in $\Theta$ and $L_n$ a positive number such that $L_n \to \infty$ as $n \to \infty$. The details of the latter approach are available at the web site given at the beginning of Section 5.

3.1. *ANOVA model.* We first spell out $\mathrm{E}(h'A^{-1}a)^2$ for the ANOVA model with normality of $v$ and $e$. We assume that the elements of the coefficient vectors $l$ and $m$ defining $\mu = l'\beta + m'v$ are bounded, and $|m| = O(1)$. In fact, $m$ typically consists of only a finite number of elements equal to 1 and the rest equal to 0. For the balanced case $\tilde{\beta} = (X'X)^{-1}X'y$ does not depend on $\Sigma$. In this case $h(\sigma) = \partial t(\sigma)/\partial \sigma = [\nabla s(\sigma)]'(I - \tilde{P}_X)u = [\nabla s(\sigma)]'u +$ lower order terms, where $u = Zv + e$, $\nabla s(\sigma) = \partial s(\sigma)/\partial \sigma'$ and $\tilde{P}_X = X(X'X)^{-1}X'$. For the general unbalanced case $\tilde{\beta}$ depends on $\Sigma$ but the same expression still holds, that is, $\partial t(\sigma)/\partial\sigma = [\nabla s(\sigma)]'u +$ lower order terms. Using Lemma 3.1 below on higher moments of normal variables, we get

$$(3.4) \quad \begin{aligned} \mathrm{E}(h'A^{-1}a)^2 &= \mathrm{tr}\left\{[\nabla s(\sigma)]'\Sigma[\nabla s(\sigma)]A^{-1}\right\} + o(d_*^{-2}) \\ &= g_3(\sigma) + o(d_*^{-2}), \end{aligned}$$

when $a$ is taken as $\partial l_\mathrm{R}(\sigma)/\partial\sigma$ corresponding to REML. Note that $d_i^2$ is the diagonal element of the information matrix associated with $\hat{\sigma}_i$ and represents the "effective sample size" for estimating $\sigma_i$.

LEMMA 3.1. *Let $A_1$ and $A_2$ be $k \times k$ symmetric matrices, and $u \sim N(0, \Sigma)$, where $\Sigma$ is $k \times k$ and positive definite. Then the following hold:*

(i) $\mathrm{E}[u\{u'A_ju - \mathrm{E}(u'A_ju)\}u'] = 2\Sigma A_j\Sigma$, $j = 1, 2$.

(ii) $\mathrm{E}[\{u'A_1u - \mathrm{E}(u'A_1u)\}\{u'A_2u - \mathrm{E}(u'A_2u)\}] = 2\,\mathrm{tr}(A_1\Sigma A_2\Sigma)$.

(iii) $\mathrm{E}[u\{u'A_1u - \mathrm{E}(u'A_1u)\}\{u'A_2u - \mathrm{E}(u'A_2u)\}u'] = 2\,\mathrm{tr}(A_1\Sigma A_2\Sigma)\Sigma + 4\Sigma A_1\Sigma A_2\Sigma + 4\Sigma A_2\Sigma A_1\Sigma$.

(iv) *Write $w_j = \lambda_j'u$, $W_j = u'A_ju$, $j = 1, \ldots, s$, where $\lambda_j$ and $A_j$ are nonstochastic of order $k \times 1$ and $k \times k$, respectively. Then, for $w = (w_1, \ldots, w_s)'$ and $W = (W_1, \ldots, W_s)'$,*

$$\begin{aligned} &\mathrm{E}[w'(W - \mathrm{E}W)]^2 \\ &= \mathrm{tr}(\Sigma_w\Sigma_W) + 4\sum_{j=1}^{s}\sum_{l=1}^{s}\lambda_j'\Sigma(A_j\Sigma A_l + A_l\Sigma A_j)\Sigma\lambda_l, \end{aligned}$$

*where $\Sigma_w$ and $\Sigma_W$ denote the covariance matrices of $w$ and $W$, respectively.*

The proofs of (i)–(iv) are immediate from Lemmas A.1 and A.2 of Prasad and Rao (1990).

It can be shown that (3.4) is also valid when $a$ is taken as $\partial l_\mathrm{P}(\sigma)/\partial\sigma$ corresponding to ML. Thus, using (1.4), (3.1), (3.3) and (3.4), we get

$$(3.5) \quad \mathrm{MSE}[t(\hat{\sigma})] = g_1(\sigma) + g_2(\sigma) + g_3(\sigma) + o(d_*^{-2}),$$

valid for both the REML estimator and the ML estimator of $\sigma$.

Next, we show that the key condition (3.2) on $t(\sigma)$ is satisfied for all balanced ANOVA models. Note that, based on the expression given by Proposition 3.1 below, all the other conditions of Theorem 3.1 are either trivial or known to be satisfied in the balanced case with normality [see Jiang (1996)]. Of course, the verification of (3.2) does not require normality.

A balanced $w$-factor linear mixed ANOVA model may be expressed as

$$(3.6) \quad y = X\beta + \sum_{i\in S}Z_iv_i + e,$$

where $X$ and $Z_i$'s have the following expressions [e.g., Searle, Casella and McCulloch (1992), Rao and Kleffe (1988)]: $X = \otimes_{l=1}^{w+1} \mathbf{1}_{n_l}^{s_l}$ with $(s_1, \ldots, s_{w+1}) \in S_{w+1} = \{0,1\}^{w+1}$, $Z_i = \otimes_{l=1}^{w+1} \mathbf{1}_{n_l}^{i_l}$ with $(i_1, \ldots, i_{w+1}) \in S \subset S_{w+1}$, where $\otimes$ denotes the operation of a Kronecker product, $n_l$ is the number of levels for factor $l$, $\mathbf{1}_n$ represents the $n$-dimensional vector of 1's, $\mathbf{1}_n^0 = I_n$ and $\mathbf{1}_n^1 = \mathbf{1}_n$. The $(w+1)$st factor corresponds to "repetition within cells," and thus $s_{w+1} = 1$ and $i_{w+1} = 1, i \in S$. We use 0 for the element $(0, \ldots, 0)$ in $S_{w+1}$ and let $\bar{S} = \{0\} \cup S$. The covariance matrix of $y$ then has the form

$$
\begin{aligned}
(3.7) \quad \Sigma &= \sigma_0 I_n + \sum_{i \in S} \sigma_i Z_i Z_i' \\
&= \sum_{i \in S_{w+1}} \lambda_i \bigotimes_{l=1}^{w+1} J_{n_l}^{i_l},
\end{aligned}
$$

where $J_n$ represents the $n \times n$ matrix of 1's, $J_n^0 = I_n$ and $J_n^1 = J_n$; $\lambda_i = \sigma_i$ if $i \in \bar{S}$, and $\lambda_i = 0$ if $i \notin \bar{S}$.

Searle and Henderson (1979) showed that $\Sigma^{-1}$ has the same form,

$$
(3.8) \qquad \Sigma^{-1} = \sum_{i \in S_{w+1}} \tau_i \bigotimes_{l=1}^{w+1} J_{n_l}^{i_l},
$$

where the coefficients $\tau_i$ in (3.8) are computed by an algorithm. From a computational point of view, the Searle–Henderson algorithm is easy to operate. However, with such a form it may not be so easy to investigate theoretical properties of $\Sigma^{-1}$, which are important to the current paper. Jiang (2004) gives an alternative derivation of (3.8) with explicit expressions for the $\tau_i$'s; see Lemma 3.2. First, note that under the balanced model we have $r_i = \prod_{i_l=0} n_l, i \in S$. This allows us to extend the definition of $r_i$ to all $i \in S_{w+1}$. In particular, $p = r_s = \prod_{s_l=0} n_l$, and $n = r_0 = \prod_{l=1}^{w+1} n_l$. We shall use the abbreviations $\{k_l = 1\}$ and so on for subsets of $L = \{1, \ldots, w+1\}$. For example, if $k, u \in S_{w+1}$, then $\{k_l = 1, u_l = 0\}$ means $\{l \in L : k_l = 1, u_l = 0\}$. Also, for $i, j \in S_{w+1}$, $j \leq i$ means $j_l \leq i_l$, $1 \leq l \leq w+1$. Finally, $|B|$ denotes the cardinality of set $B$.

LEMMA 3.2.   *For any balanced $w$-factor mixed ANOVA model* (3.8) *holds with*

$$
(3.9) \quad \tau_i = \left(\frac{r_i}{n}\right) \left\{ \sum_{j \in S_{w+1}} \frac{(-1)^{|\{i_l = 1, j_l = 0\}|} \mathbb{1}_{(j \leq i)}}{\sigma_0 + \sum_{k \in S} \sigma_k (n/r_k) \mathbb{1}_{(k \leq j)}} \right\}, \qquad i \in S_{w+1}.
$$

Using Lemma 3.2, the following proposition can be proved. A sketch of the proof is given in Section 5.5.

PROPOSITION 3.1. *For any balanced mixed ANOVA model, the BLUP* $t(\sigma)$ *given by* (1.3) *can be expressed as* (3.2) *with* $K \leq 1 + |S|2^{w+1}$ *[hence* $K = O(1)$ *and the terms below* (3.2) *bounded].*

It is known that (3.2) also holds for some unbalanced ANOVA models. For example, see the nested error regression model discussed in the next subsection, which is also a special case of the longitudinal model.

3.2. *The longitudinal model.* For longitudinal models the spelled-out expression for $\mathrm{E}(h'A^{-1}a)^2$ in (3.3) [up to a term $o(d_*^{-2})$] is given in Datta and Lahiri (2000). In the following we show that the key condition (3.2) in Theorem 3.1 is satisfied for two special (and important) classes of longitudinal models: the Fay–Herriot model and the nested error regression model.

First consider the Fay–Herriot model (see Section 2.3). It is easy to show the following:

$$(3.10) \qquad l'\tilde{\beta} = \sum_{i=1}^{t} a_i(\sigma)y_i,$$

$$m'\tilde{v} = \sum_{i=1}^{t} m_i \tilde{v}_i$$

$$(3.11) \qquad = \sum_{i=1}^{t} b_i(\sigma)y_i - \sum_{i,k=1}^{t} b_{i,k}(\sigma)y_k,$$

where

$$a_i(\sigma) = l'\left(\sum_{i=1}^{t} \frac{x_i x_i'}{\sigma + \phi_i}\right)^{-1} \frac{x_i}{\sigma + \phi_i},$$

$$b_i(\sigma) = m_i\left(\frac{\sigma}{\sigma + \phi_i}\right),$$

$$b_{i,k} = m_i\left(\frac{\sigma}{\sigma + \phi_i}\right)x_i'\left(\sum_{k=1}^{t} \frac{x_k x_k'}{\sigma + \phi_k}\right)^{-1}\left(\frac{\sigma}{\sigma + \phi_k}\right).$$

We have the following result. The proof is straightforward.

PROPOSITION 3.2. *For the Fay–Herriot model* (2.14), *the BLUP* $t(\sigma)$ *given by* (1.3) *can be expressed as* (3.2) *with* $K = (t+1)^2 - 1$ *[and the terms below* (3.2) *bounded], provided that* (i) *the* $\phi_i$'s *are bounded from above and away from* 0; (ii) $|x_i|$, $1 \leq i \leq t$, *are bounded, and so are* $|l|$ *and* $\sum_{i=1}^{t} |m_i|$, *and* (iii) $\liminf \lambda_{\min}(t^{-1}\sum_{i=1}^{t} x_i x_i') > 0$.

Next we consider the nested error regression model. Suppose that

$$(3.12) \qquad y_{ij} = \beta_0 + x_{ij}'\beta + v_i + e_{ij}, \qquad j = 1, \ldots, n_i;\ i = 1, \ldots, t,$$

where $\beta = (\beta_u)_{1 \le u \le p-1}$ and $\beta_u$, $0 \le u \le p-1$, are unknown regression coefficients; $x_{ij}$'s are known vectors of covariates; $v_i$'s are independent random effects with distribution $N(0, \sigma_1)$; $e_{ij}$'s are independent errors with distribution $N(0, \sigma_0)$, and $v$ and $e$ are independent. W.l.o.g. let $n_i \ge 1$. For the nested error regression model (3.12) $l'\tilde{\beta}$ and $m'\tilde{v}$ can, again, be expressed as (3.10) and (3.11), where

$$a_i(\sigma) = l' \left( \sum_{i=1}^{t} X_i' \Sigma_i^{-1} X_i \right)^{-1} X_i' \Sigma_i^{-1},$$

$$b_i(\sigma) = m_i \left( \frac{\sigma_1}{\lambda_i} \right) \mathbf{1}_{n_i}',$$

$$b_{i,k}(\sigma) = m_i \left( \frac{\sigma_1}{\lambda_i} \right) \mathbf{1}_{n_i}' X_i \left( \sum_{k=1}^{t} X_k' \Sigma_k^{-1} X_k \right)^{-1} X_k' \Sigma_k^{-1}.$$

We have the following result. A sketch of the proof is given in Section 5.7.

PROPOSITION 3.3.   *For the nested error regression model* (3.12) *the BLUP* $t(\sigma)$ *given by* (1.3) *can be expressed as* (3.2) *with* $K = O(t^2)$ [*and the terms below* (3.2) *bounded*], *provided that* (i) $\sigma_0 > 0$; (ii) $p$, $n_i$, $|x_{ij}|$, $1 \le i \le t$, $1 \le j \le n_i$, *are bounded, and so are* $|l|$ *and* $\sum_{i=1}^{t} |m_i|$, *and* (iii) $\liminf \lambda_{\min}(t^{-a} S_a) > 0$, $a = 1, 2$, *where* $S_1 = \sum_{i=1}^{t} n_i \sum_{j=1}^{n_i} (x_{ij} - \bar{x}_i)(x_{ij} - \bar{x}_i)'$, $S_2 = \sum_{i=1}^{t} \sum_{j=1}^{n_i} \sum_{k=1}^{t} \sum_{l=1}^{n_k} (x_{ij} - x_{kl})(x_{ij} - x_{kl})'$ *and* $\bar{x}_i = n_i^{-1} \sum_{j=1}^{n_i} x_{ij}$.

Note that in both cases considered above $d_*$ can be chosen as $\sqrt{t}$.

**4. Estimation of MSE.**   We now turn to the estimation of $\text{MSE}[t(\hat{\sigma})]$. We obtain an estimator $\text{mse}[t(\hat{\sigma})]$ correct to second order in the sense of $E\{\text{mse}[t(\hat{\sigma})]\} = \text{MSE}[t(\hat{\sigma})] + o(d_*^{-2})$. That is, the bias of $\text{mse}[t(\hat{\sigma})]$ in estimating $\text{MSE}[t(\hat{\sigma})]$ is $o(d_*^{-2})$.

First, we have from (3.1) and (3.3),

$$\begin{aligned}
(4.1) \qquad \text{MSE}[t(\hat{\sigma})] &= \text{MSE}[t(\sigma)] + E(h'A^{-1}a)^2 + o(d_*^{-2}) \\
&= \eta(\sigma) + o(d_*^{-2}), \qquad \text{say.}
\end{aligned}$$

We now define an estimator $\hat{\eta}$ of $\eta(\sigma)$ having the following property:

$$(4.2) \qquad\qquad E(\hat{\eta}) = \eta(\sigma) + o(d_*^{-2}).$$

It follows from (4.1) and (4.2) that the bias of $\hat{\eta}$ in estimating $\text{MSE}[t(\hat{\sigma})]$ is $o(d_*^{-2})$. In addition to $a$ and $A$ defined in Section 2 (Theorem 2.1), let $b = \partial \eta(\sigma) / \partial \sigma = (b_i)$; $B = \partial^2 \eta(\sigma) / \partial \sigma^2 = (b_{ij})$; $F = \partial^2 l(\sigma) / \partial \sigma^2$, $H_i = (\partial^3 l(\sigma) / \partial \sigma_i \partial \sigma^2)$ and $C = (a'A^{-1} h_i)_{1 \le i \le s}$, where $s$ is the dimension of $\sigma$.

Also, let $Q = D^{-1}AD^{-1}$ and $W = Q^{-1} = (w_{ij})$. Let $D^{-1}a = (\lambda_i)$, $D^{-1/2}(F - A)D^{-1/2} = (\lambda_{ij})$ and $D^{-1}H_iD^{-1} = (\lambda_{ijk})$. We define the following vector, matrix and arrays: $U_0 = (u_i)$, $U_1 = (u_{il})$, $U_2 = (u_{jkl})$ and $U_3 = (u_{jklmn})$, where $u_i = \mathrm{E}(\lambda_i)$, $u_{il} = \mathrm{E}(\lambda_i\lambda_l)$, $u_{jkl} = \mathrm{E}(\lambda_{jk}\lambda_l)$ and $u_{jklmn} = \mathrm{E}(\lambda_{jkm}\lambda_l\lambda_n)$. Note that all of these are functions of $\sigma$ [e.g., $A = A(\sigma)$]. The norm of an $r$-way array ($r \geq 3$) $U$, denoted by $\|U\|$, is defined as the maximum of the absolute values of its elements. (The norm of a matrix is defined in Theorem 3.1.) Define

$$(4.3) \qquad \Delta_0(\sigma) = -2b'A^{-1}\mathrm{E}(a),$$

$$(4.4) \qquad \Delta_1(\sigma) = b'A^{-1}\mathrm{E}(FA^{-1}a),$$

$$(4.5) \qquad \Delta_2(\sigma) = \tfrac{1}{2}\mathrm{E}(a'A^{-1}BA^{-1}a),$$

$$(4.6) \qquad \Delta_3(\sigma) = -\tfrac{1}{2}b'A^{-1}\mathrm{E}(CA^{-1}a).$$

Finally, we define

$$(4.7) \qquad \hat{\eta} = \eta(\hat{\sigma}) - \sum_{j=0}^{3}\Delta_j(\hat{\sigma}),$$

provided that $|\hat{\eta}| \leq c_0 d_*^\lambda$; otherwise, let $\hat{\eta} = \eta(\sigma^*)$, where $c_0$ and $\lambda$ are known positive constants and $\sigma^*$ is a given point in $\Theta$.

THEOREM 4.1. *The estimator $\hat{\eta}$ given above satisfies the property* (4.2) *provided that*

(i) *$\eta(\cdot)$ is three times continuously differentiable and the following are bounded: $\eta(\sigma)$, $|b|$, $\|B\|$ and*

$$\sup_{\tilde{\sigma} \in S_\delta(\sigma)} \left| \frac{\partial^3 \eta(\tilde{\sigma})}{\partial\sigma_i\,\partial\sigma_j\,\partial\sigma_k} \right|, \qquad 1 \leq i,\ j, k \leq s,$$

*where $\delta$ is a positive number and $S_\delta(\sigma) = \{\tilde{\sigma} : |\tilde{\sigma}_i - \sigma_i| \leq \delta d_*/d_i, 1 \leq i \leq s\}$.*

(ii) *The conditions of Theorem 2.1 hold with $g > 8 + 4\lambda$ and $l(\sigma)$ four times continuously differentiable with respect to $\sigma$.*

(iii) *The $g$th moments of the following are bounded:*

$$\frac{1}{\sqrt{d_j d_k}} \left| \frac{\partial^3 l(\sigma)}{\partial\sigma_i\,\partial\sigma_j\,\partial\sigma_k} - \mathrm{E}\left[ \frac{\partial^3 l(\sigma)}{\partial\sigma_i\,\partial\sigma_j\,\partial\sigma_k} \right] \right|, \qquad \frac{1}{d_j d_k} \left| \frac{\partial^3 l(\sigma)}{\partial\sigma_i\,\partial\sigma_j\,\partial\sigma_k} \right|$$

*and*

$$\frac{d_*^2}{d_i d_j d_k d_l} \sup_{\tilde{\sigma} \in S_\delta(\sigma)} \left| \frac{\partial^4 l(\tilde{\sigma})}{\partial\sigma_i\,\partial\sigma_j\,\partial\sigma_k\,\partial\sigma_l} \right|, \qquad 1 \leq i,\ j, k, l \leq s.$$

(iv) *$\sup_{\tilde{\sigma} \in S_\delta(\sigma)} \|Q(\tilde{\sigma}) - Q(\sigma)\| \to 0$ and $\sup_{\tilde{\sigma} \in S_\delta(\sigma)} \|U_j(\tilde{\sigma}) - U_j(\sigma)\| \to 0$, $j = 1, 2, 3$, as $\delta \to 0$ uniformly in $n$.*

(v) $|\mathrm{E}(a)|$ *is bounded and* $\sup_{\tilde{\sigma}\in S_\delta(\sigma)}|\mathrm{E}(a)|_{\sigma=\tilde{\sigma}}-\mathrm{E}(a)|\to 0$ *as* $\delta\to 0$ *uniformly in* $n$.

A sketch of the proof is given in Section 5.8.

COROLLARY 4.1.  *If condition* (v) *of Theorem* 4.1 *is strengthened to*

$$(4.8)\qquad\qquad\qquad \mathrm{E}(a)=0,$$

*then* (4.2) *holds, in which* $\hat{\eta}$ *is given as in Theorem* 4.1 *except that in* (4.7) *the summation is from* 1 *to* 3.

The estimator $\hat{\eta}$ considered in Theorem 4.1 is truncated when its value exceeds some (large) threshold. Such a truncation is needed only for establishing the asymptotic result. In practice one does not need to truncate the estimator (because it can be argued that for a given value of $\hat{\eta}$ there are always constants $c_0$ and $\lambda$ such that $|\hat{\eta}|\le c_0 d_*^\lambda$). On the other hand, a similar result may be obtained for $\hat{\eta}$ without truncation, provided that $\hat{\sigma}$ is replaced by its truncated version (defined above Section 3.1). The details of such a result are available at the web site given at the beginning of Section 5.

4.1. *REML and ML under the ANOVA model.*  In Section 5.9 we give sketches of a proof that shows the conditions of Theorem 4.1 are satisfied if the REML estimator of $\sigma$ is used, provided that the same conditions under which the REML estimators are consistent [Jiang (1996)] are satisfied. It can be shown that similar results also hold for ML estimation, but we omit the details. According to (3.5), in both REML and ML cases we have $\eta(\sigma)=g_1(\sigma)+g_2(\sigma)+g_3(\sigma)$. However, there is a difference between the two in terms of the spelled-out formulas for MSE estimation. This difference is due to the fact that (4.8) holds for REML but not for ML.

First consider REML. Letting $a=a_\mathrm{R}$, $l(\sigma)=l_\mathrm{R}(\sigma)$ and $A=A_\mathrm{R}$, we have $\mathrm{E}(a_\mathrm{R})=0$, and the $(i,j)$ element of $A_\mathrm{R}$ is given by $-(1/2)\mathrm{tr}(PV_iPV_j)$. Furthermore, we have $\Delta_0(\sigma)=0$, $\Delta_1(\sigma)=b'A_\mathrm{R}^{-1}w_\mathrm{R}(\sigma)+o(d_*^{-2})$, where $w_\mathrm{R}(\sigma)=(w_{0,\mathrm{R}},\ldots,w_{q,\mathrm{R}})'$ with $w_{i,\mathrm{R}}=-\mathrm{tr}\{A_\mathrm{R}^{-1}[\mathrm{tr}(PV_iPV_jPV_k)]_{0\le j,k\le q}\}$, $i=0,\ldots,q$; $\Delta_2(\sigma)=-g_3(\sigma)+o(d_*^{-2})$, where $g_3(\sigma)$ is given by (3.4) with $A=A_\mathrm{R}$ and $\Delta_3(\sigma)=-b'A_\mathrm{R}^{-1}w_\mathrm{R}(\sigma)+o(d_*^{-2})$. Thus for REML we have $\sum_{j=0}^3\Delta_j(\sigma)=-g_3(\sigma)+o(d_*^{-2})$. It follows from (4.7) that for REML $\hat{\eta}=\hat{\eta}_\mathrm{R}$, where

$$(4.9)\qquad\qquad \hat{\eta}_\mathrm{R}=g_1(\hat{\sigma}_\mathrm{R})+g_2(\hat{\sigma}_\mathrm{R})+2g_3(\hat{\sigma}_\mathrm{R}),$$

where $\hat{\sigma}_\mathrm{R}$ is the REML estimator of $\sigma$. The MSE estimator $\hat{\eta}_\mathrm{R}$ given by (4.9) depends on the data $y$ only through $\hat{\sigma}_\mathrm{R}$. An alternative MSE estimator that is data specific can be obtained by using the following estimator of $g_3(\sigma)$:

$$(4.10)\qquad \tilde{g}_3(\hat{\sigma}_\mathrm{R})=(y-X\tilde{\beta})'[\nabla s(\sigma)]'A_\mathrm{R}^{-1}[\nabla s(\sigma)](y-X\tilde{\beta})|_{\sigma=\hat{\sigma}_\mathrm{R}},$$

where $\tilde{\beta}$ is the BLUE given below (1.3). It can be shown that $\mathrm{E}[\tilde{g}_3(\hat{\sigma}_{\mathrm{R}})] = g_3(\sigma) + o(d_*^{-2})$. The estimator (4.10) is obtained from (3.4) by replacing $\Sigma$ by $\hat{\Sigma} = (y - X\hat{\beta}_{\mathrm{R}})(y - X\hat{\beta}_{\mathrm{R}})'$.

Now consider ML. For simplicity, we assume that $\mathrm{rank}(X) = p$ is bounded. Then, similarly, letting $a = a_{\mathrm{M}}$, $l(\sigma) = l_{\mathrm{P}}(\sigma)$ and $A = A_{\mathrm{M}}$ for ML, we have $\mathrm{E}(a_{\mathrm{M}}) = -g_{\mathrm{M}0}(\sigma) = -g_{\mathrm{M}0}$ with $g_{\mathrm{M}0,i} = (1/2)\,\mathrm{tr}[(\Sigma^{-1} - P)V_i]$, $i = 0, \ldots, q$. Furthermore, $\Delta_0(\sigma) = 2b'A_{\mathrm{M}}^{-1}g_{\mathrm{M}0}$; $\Delta_1(\sigma) = b'A_{\mathrm{M}}^{-1}w_{\mathrm{M}} - b'A_{\mathrm{M}}^{-1}g_{\mathrm{M}0} + o(d_*^{-2})$, where $w_{\mathrm{M}}$ is $w_{\mathrm{R}}$ with $P$ replaced by $\Sigma^{-1}$; $\Delta_2(\sigma) = -g_3(\sigma) + o(d_*^{-2})$ and $\Delta_3(\sigma) = -b'A_{\mathrm{M}}^{-1}w_{\mathrm{M}} + o(d_*^{-2})$. Also, $b'A_{\mathrm{M}}^{-1}g_{\mathrm{M}0} = (\partial g_1/\partial\sigma)'A_{\mathrm{M}}^{-1}g_{\mathrm{M}0} + o(d_*^{-2}) = g_{10}(\sigma) + o(d_*^{-2})$, say. Thus for ML we have $\sum_{j=0}^{3}\Delta_j(\sigma) = g_{10}(\sigma) - g_3(\sigma) + o(d_*^{-2})$, hence $\hat{\eta} = \hat{\eta}_{\mathrm{M}}$, where

$$(4.11) \qquad \hat{\eta}_{\mathrm{M}} = g_1(\hat{\sigma}_{\mathrm{M}}) + g_2(\hat{\sigma}_{\mathrm{M}}) + 2g_3(\hat{\sigma}_{\mathrm{M}}) - g_{10}(\hat{\sigma}_{\mathrm{M}}),$$

where $\hat{\sigma}_{\mathrm{M}}$ is the ML estimator of $\sigma$. Similar to the REML case, a data-specific estimator can be obtained by using $\tilde{g}_3(\hat{\sigma}_{\mathrm{M}})$ instead of $g_3(\hat{\sigma}_{\mathrm{M}})$.

4.2. *REML and ML under the longitudinal model.* For longitudinal models the spell-out of (4.7) [up to a term $o(d_*^{-2})$] is given by Datta and Lahiri (2000) for REML and ML estimation. Note that, similar to the previous subsection, there is a difference between using REML and ML. In the following we give regularity conditions such that the conditions of Theorem 4.1 are satisfied for longitudinal models. The assumption that $G_i$ and $R_i$ are linear functions of $\sigma$ can be relaxed.

We consider REML estimation. Similar results also hold for ML but we shall omit the details. Let $\sigma = (\sigma_0, \ldots, \sigma_q)'$. Suppose that

1. $G_i$ and $R_i$ are linear in $\sigma$ such that $\|G_i\|$, $\|R_i\|$, $\|\partial G_i/\partial\sigma_j\|$, $\|\partial R_i/\partial\sigma_j\|$, $0 \le j \le q$, and $\|\Sigma_i^{-1}\|$ are bounded, and, as $\tilde{\sigma} \to \sigma$, $\max_{1 \le i \le t}\|G_i(\tilde{\sigma}) - G_i(\sigma)\| \to 0$, $\|R_i(\tilde{\sigma}) - R_i(\sigma)\| \to 0$ uniformly in $t$:
2. $\sigma \in \Theta^{\mathrm{o}}$, the interior of $\Theta$.
3. The following are bounded: $p$, $n_i$, $\|X_i\|$, $\|Z_i\|$, $|l|$ and $\sum_{i=1}^{t}|m_i|$.
4. $\liminf_{t\to\infty}\lambda_{\min}(t^{-1}\sum_{i=1}^{t}B_i) > 0$, $\liminf_{t\to\infty}\lambda_{\min}(t^{-1}\sum_{i=1}^{t}X_i'\Sigma_i^{-1}X_i) > 0$, where $B_i$ is the $(q+1)\times(q+1)$ matrix whose $(j,k)$ element is $\mathrm{tr}(\Sigma_i^{-1}\Sigma_{i,j} \times \Sigma_i^{-1}\Sigma_{i,k})$ with $\Sigma_{i,j} = \partial\Sigma_i/\partial\sigma_j$. Then the conditions of Theorem 4.1 are satisfied for the longitudinal model. A sketch of the proof is given in Section 5.10.

## 5. Sketches of proofs.

In this section we give very brief sketches of the proofs involved in this paper. These include proofs of the theorems and other technical details. The detailed proofs are available at the following web site address: http://anson.ucdavis.edu/˜jiang/jp8.pdf.

5.1. *Regarding Theorem* 2.1. Let $\sigma_* = \sigma - A^{-1}a$ and $\mathcal{B} = B_1 \cap B_2$, where $B_1 = \{|\xi| \leq (|\lambda|/2)d_*^{1-\rho}\}$ and $B_2 = \{\|\eta\|(1+|\lambda|^{-1}|\xi|)^2 + (s^{3/2}/3)\zeta(1+|\lambda|^{-1}|\xi|)^3 < |\lambda|/2\}$ with $\xi = D^{-1}a$, $\lambda = \lambda_{\max}(D^{-1}AD^{-1})$, $\eta = D^{-1}(\partial^2 l(\sigma)/\partial\sigma^2 - A)D^{-1}$ and $\zeta = \max_{i,j,k}\{M_{ijk}(y)/d_i d_j d_k\}$. It can be shown that, on $\mathcal{B}$, $l(\tilde{\sigma}) - l(\sigma_*) < 0$ if $|D(\tilde{\sigma} - \sigma_*)| = 1$. Since the function $l(\tilde{\sigma})$ cannot attain its maximum over the set $\{\tilde{\sigma} : |D(\tilde{\sigma} - \sigma_*)| \leq 1\}$ at the boundary of the set, the maximum must be attained in the interior. Thus, there is a solution to (2.1) in $\{\tilde{\sigma} : |D(\tilde{\sigma} - \sigma_*)| < 1\}$. Let $\hat{\sigma}$ be the solution to (2.1) closest to $\sigma_*$. It can be shown that, on $\mathcal{B}$, $\hat{\sigma} \in \Theta$, $\partial l(\hat{\sigma})/\partial\sigma = 0$ and $|D(\hat{\sigma} - \sigma)| < d_*^{1-\rho}$. Furthermore, by Taylor expansion it can be shown that, on $\mathcal{B}$, $\hat{\sigma} - \sigma = -A^{-1}a + r_1 + r_2$ such that $|r_1| \leq d_*^{-1-\rho}u_1$, $|r_2| \leq d_*^{-2\rho}u_2$ with $\mathrm{E}(u_j^g)$, $j = 1, 2$, bounded. Finally, it can be shown that $\mathrm{P}(B_1^c) \leq c_1 d_*^{-(1-\rho)g}$ and $\mathrm{P}(B_2^c) \leq c_2 d_*^{-g/4}$ for some constants $c_1$ and $c_2$.

5.2. *Regarding Section* 2.1. The following lemmas are used to verify that the conditions of Theorem 2.1 are satisfied.

LEMMA 5.1. *Let $Q$ be a symmetric matrix and $\xi \sim N(0, I)$. Then for any $g \geq 2$ there is a constant $c$ that only depends on $g$ such that $\mathrm{E}|\xi'Q\xi - \mathrm{E}\xi'Q\xi|^g \leq c\|Q\|_2^g$.*

LEMMA 5.2. *For any matrices $A$, $B$ and $C$, we have $|\mathrm{tr}(ABC)| \leq \|B\| \cdot \|A\|_2 \cdot \|C\|_2$, provided that the matrix product is well defined.*

LEMMA 5.3. *Let $a_i$, $b_i$ be real numbers and $\delta_i \geq 0$ such that $a_i \leq b_i + \Delta_a$, $1 \leq i \leq s$, where $\Delta_a = \sum_{j=1}^s \delta_j a_j$. If $\Delta = \sum_{j=1}^s \delta_j < 1$, then $a_i \leq b_i + (1-\Delta)^{-1}\Delta_b$, $1 \leq i \leq s$, where $\Delta_b = \sum_{j=1}^s \delta_j b_j$.*

We also use the following expressions for second and third derivatives of $l_{\mathrm{R}}(\sigma)$:

$$(5.1) \qquad \frac{\partial^2 l_{\mathrm{R}}(\sigma)}{\partial\sigma_i\,\partial\sigma_j} = \frac{1}{2}\,\mathrm{tr}(PV_iPV_j) - u'PV_iPV_jPu,$$

$$\frac{\partial^3 l_{\mathrm{R}}(\sigma)}{\partial\sigma_i\,\partial\sigma_j\,\partial\sigma_k} = u'PV_iPV_jPV_kPu$$

$$+ u'PV_jPV_kPV_iPu + u'PV_kPV_iPV_jPu$$

$$(5.2) \qquad - \frac{1}{2}[\mathrm{tr}(PV_iPV_jPV_k) + \mathrm{tr}(PV_iPV_kPV_j)].$$

Note that the second and third derivatives are involved in the conditions of Theorem 2.1.

5.3. *Regarding Section* 2.2. In addition to (2.13), we have

$$(5.3) \qquad \frac{\partial^2 l_{\mathrm{P}}(\sigma)}{\partial \sigma_i \partial \sigma_j} = \frac{1}{2} \operatorname{tr}(\Sigma^{-1} V_i \Sigma^{-1} V_j) - u' P V_i P V_j P u,$$

$$\frac{\partial^3 l_{\mathrm{P}}(\sigma)}{\partial \sigma_i \partial \sigma_j \partial \sigma_k} = u' P V_i P V_j P V_k P u$$

$$+ u' P V_j P V_k P V_i P u + u' P V_k P V_i P V_j P u$$

$$(5.4) \qquad\qquad\qquad - \frac{1}{2} [\operatorname{tr}(\Sigma^{-1} V_i \Sigma^{-1} V_j \Sigma^{-1} V_k)$$

$$+ \operatorname{tr}(\Sigma^{-1} V_i \Sigma^{-1} V_k \Sigma^{-1} V_j)].$$

From these expressions we obtain the following relationships:

$$\frac{\partial l_{\mathrm{P}}(\sigma)}{\partial \sigma_i} = \frac{\partial l_{\mathrm{R}}(\sigma)}{\partial \sigma_1} + \operatorname{tr}(P V_i) - \operatorname{tr}(\Sigma^{-1} V_i),$$

$$(5.5) \qquad \frac{\partial^2 l_{\mathrm{P}}(\sigma)}{\partial \sigma_i \partial \sigma_j} = \frac{\partial^2 l_{\mathrm{R}}(\sigma)}{\partial \sigma_i \partial \sigma_j} + \frac{1}{2} [\operatorname{tr}(\Sigma^{-1} V_i \Sigma^{-1} V_j) - \operatorname{tr}(P V_i P V_j)],$$

$$\frac{\partial^3 l_{\mathrm{P}}(\sigma)}{\partial \sigma_i \partial \sigma_j \partial \sigma_k} = \frac{\partial^3 l_{\mathrm{R}}(\sigma)}{\partial \sigma_i \partial \sigma_j \partial \sigma_k} + \frac{1}{2}(I_1 - J_1 + I_2 - J_2),$$

where $I_1 = \operatorname{tr}(P V_i P V_j P V_k)$, $I_2 = \operatorname{tr}(P V_i P V_k P V_j)$ and $J_r$ is $I_r$ with $P$ replaced by $\Sigma^{-1}$, $r = 1, 2$. We assume that $p = \operatorname{rank}(X)$ is bounded. Then it can be shown that $|\operatorname{tr}(\Sigma^{-1} V_i) - \operatorname{tr}(P V_i)| \le p\sigma_i^{-1}$, $|\operatorname{tr}(\Sigma^{-1} V_i \Sigma^{-1} V_j) - \operatorname{tr}(P V_i P V_j)| \le 3p(\sigma_i \sigma_j)^{-1}$ and $|\operatorname{tr}(\Sigma^{-1} V_i \Sigma^{-1} V_j \Sigma^{-1} V_k) - \operatorname{tr}(P V_i P V_j P V_k)| \le 7p(\sigma_i \sigma_j \sigma_k)^{-1}$. Thus, by the result of the previous subsection, it can be shown that the conditions of Theorem 2.1 are satisfied.

5.4. *Regarding Theorem* 3.1. Let $\rho = 3/4$. By Theorem 2.1 and Taylor expansion it can be shown that $t(\hat{\sigma}) - t(\sigma) = -h' A^{-1} a + r$, where $|r| \le d_*^{-2\rho} u$ and $\mathrm{E}(u^2)$ is bounded. Thus, we have $\mathrm{E}[t(\hat{\sigma}) - t(\sigma)]^2 = \mathrm{E}(\cdot)^2 \mathbf{1}_{\mathcal{B}} + \mathrm{E}(\cdot)^2 \mathbf{1}_{\mathcal{B}^c}$, where $(\cdot)^2$ denotes $[t(\hat{\sigma}) - t(\sigma)]^2$. The first term $= \mathrm{E}(h' A^{-1} a)^2 \mathbf{1}_{\mathcal{B}} + O(d_*^{-(1+2\rho)}) + O(d_*^{-4\rho})$, while $\mathrm{E}(h' A^{-1} a)^2 \mathbf{1}_{\mathcal{B}^c} = O(d_*^{-2-\nu})$ for some $\nu > 0$.

5.5. *Regarding Proposition* 3.1. First, the following identity can be established:

$$(5.6) \qquad \begin{aligned} P &= \Sigma^{-1} - \Sigma^{-1} X (X' \Sigma^{-1} X)^{-1} X' \Sigma^{-1} \\ &= \left\{ I_n - \left(\frac{p}{n}\right) X X' \right\} \Sigma^{-1}. \end{aligned}$$

By the definition of BLUP for $v$, Lemma 3.2 and (5.6), it can be shown that $\tilde{v}_i = \sum_{k \in S_{w+1}} \sigma_i \tau_k W_{i,k} y$, where the $\tau_k$'s are given by (3.9) and

$$W_{i,k} = Z_i' \Big\{ I_n - \Big(\frac{p}{n}\Big) X X' \Big\} \bigotimes_{l=1}^{w+1} J_{n_l}^{k_l}, \qquad k \in S_{w+1}.$$

5.6. *Regarding the Fay–Herriot model.* For REML estimation, the restricted loglikelihood is given by $l_{\mathrm{R}}(\sigma) = c - (1/2)(\log |T' \Sigma T| + y' P y)$, where $c$ is a constant, $T$ is as in Section 2.1 and $P = T(T' \Sigma T)^{-1} T' =$ the middle term of (5.6) with $\Sigma = \sigma I + \Phi$. Suppose that $\sigma > 0$ and the $\phi_i$'s are bounded. Then it can be shown that the conditions of Theorem 2.1 are satisfied with $D = d = \sqrt{t}$. A similar result can be proved for ML estimation, in which case one considers the profile loglikelihood $l_{\mathrm{P}}(\sigma) = c - (1/2)[\log |\Sigma| + y' P y]$.

5.7. *Regarding Proposition 3.3.* First note that $X_i = (\mathbf{1}_{n_i} x_i)$, where the $j$th row of $x_i$ is $x_{ij}'$. Also, we have $\Sigma_i = \sigma_0 I_{n_i} + \sigma_1 J_{n_i}$, thus $\Sigma_i^{-1} = \lambda_i^{-1} I_{n_i} + \gamma n_i \lambda_i^{-1}(I_{n_i} - \bar{J}_{n_i})$, where $\gamma = \sigma_1/\sigma_0$, $\lambda_i = \lambda_i(\sigma) = \sigma_0 + n_i \sigma_1$ and $\bar{J}_{n_i} = J_{n_i}/n_i$. Therefore, we can write

$$\sum_{i=1}^{t} X_i' \Sigma_i^{-1} X_i = \begin{pmatrix} A & B' \\ B & C + \gamma D \end{pmatrix},$$

where $A = \sum_{i=1}^{t} n_i/\lambda_i$, $B = \sum_{i=1}^{t} x_i' \mathbf{1}_{n_i}/\lambda_i$, $C = \sum_{i=1}^{t} x_i' x_i/\lambda_i$ and $D = \sum_{i=1}^{t} (n_i/\lambda_i) x_i'(I_{n_i} - \bar{J}_{n_i}) x_i$. Thus,

$$(5.7) \qquad \Big( \sum_{i=1}^{t} X_i' \Sigma_i^{-1} X_i \Big)^{-1} = \begin{pmatrix} Q & -A^{-1} B' R \\ -A^{-1} R B & R \end{pmatrix},$$

where $Q = [A - B'(C + \gamma D)^{-1} B]^{-1}$ and $R = (C + \gamma D - A^{-1} B B')^{-1}$. It can be shown that $AC - BB' \geq S_2/2\lambda_{\mathrm{M}}^2$, where $\lambda_{\mathrm{M}} = \max_i \lambda_i = \sigma_0 + n_{\max} \sigma_1$ with $n_{\max} = \max_i n_i$. It follows, by conditions (ii) and (iii), that $\|R\| \leq \lambda_{\mathrm{M}}/t(\delta_1 \gamma + \delta_2)$, where $\delta_a$, $a = 1, 2$, are some positive constants. Then, using the identity $Q = A^{-1} + A^{-2} B' R B$, one can show $\|Q\| \leq c(\lambda_{\mathrm{M}}/t)$, where $c$ is a constant.

5.8. *Regarding Theorem 4.1.* The proof of Theorem 4.1 requires the following lemmas [see Jiang, Lahiri and Wan (2002)].

LEMMA 5.4. *For any nonsingular matrices $P$, $Q$ and nonnegative integer $q$,*

$$Q^{-1} = \Big( \sum_{r=0}^{q} [P^{-1}(P - Q)]^r \Big) P^{-1} + [P^{-1}(P - Q)]^{q+1} Q^{-1}.$$

Lemma 5.5. *Let $P$, $Q$ be matrices such that $P$ is nonsingular and $\|Q - P\| \le (3\|P^{-1}\|)^{-1}$. Then $Q$ is nonsingular and $\|Q^{-1}\| \le (3/\sqrt{2})\|P^{-1}\|$.*

Let $\mathcal{A}$ be the set such that the following hold:

$$\frac{1}{\sqrt{d_i d_j}}\left|\frac{\partial^2 l(\sigma)}{\partial \sigma_i \, \partial \sigma_j} - \mathrm{E}\left[\frac{\partial^2 l(\sigma)}{\partial \sigma_i \, \partial \sigma_j}\right]\right| \le d_*^\tau, \qquad 1 \le i, j \le s,$$

$$\frac{1}{\sqrt{d_j d_k}}\left|\frac{\partial^3 l(\sigma)}{\partial \sigma_i \, \partial \sigma_j \, \partial \sigma_k} - \mathrm{E}\left[\frac{\partial^3 l(\sigma)}{\partial \sigma_i \, \partial \sigma_j \, \partial \sigma_k}\right]\right| \le d_*^\tau, \qquad 1 \le i, j, k \le s.$$

Let $\mathcal{E} = \mathcal{A} \cap \mathcal{B}$. Let $\rho = 3/4$ in Theorem 2.1. It can be shown that $\mathrm{P}(\mathcal{E}^c) \le c d_*^{-\tau g}$, where $\tau = 1/4$. By Taylor expansion, it can be shown that the following holds on $\mathcal{E}$:

$$\eta(\hat{\sigma}) = \eta(\sigma) - 2b'A^{-1}a + b'A^{-1}fA^{-1}a$$
$$+ \tfrac{1}{2}[a'A^{-1}BA^{-1}a - b'A^{-1}CA^{-1}a] + r,$$

where $|r| \le d_*^{-3\rho}u$ and $\mathrm{E}(u)$ is bounded. Thus it can be shown that $\mathrm{E}\eta(\hat{\sigma})\mathbf{1}_\mathcal{E} = \eta(\sigma) + \sum_{j=0}^{3}\Delta_j(\sigma) + o(d_*^{-2})$. On the other hand, we have the following expressions:

$$\Delta_0(\sigma) = -2\sum_{i,j}\frac{1}{d_i}b_i(\sigma)w_{ij}(\sigma)u_j(\sigma),$$

$$\Delta_1(\sigma) = \sum_{i,j,k,l}\frac{1}{d_i\sqrt{d_j d_k}}b_i(\sigma)w_{ij}(\sigma)w_{kl}(\sigma)u_{jkl}(\sigma),$$

$$\Delta_2(\sigma) = \frac{1}{2}\sum_{i,j,k,l}\frac{1}{d_j d_k}b_{jk}(\sigma)w_{ij}(\sigma)w_{kl}(\sigma)u_{il}(\sigma),$$

$$\Delta_3(\sigma) = -\frac{1}{2}\sum_{i,j,k,l,m,n}\frac{1}{d_i d_j}b_i(\sigma)w_{ij}(\sigma)w_{kl}(\sigma)w_{mn}(\sigma)u_{jklmn}(\sigma).$$

With these it can be shown that $\mathrm{E}\Delta_j(\hat{\sigma})\mathbf{1}_\mathcal{E} = \Delta_j(\sigma) + o(d_*^{-2})$. It follows that $\mathrm{E}\hat{\eta}\mathbf{1}_\mathcal{E} = \eta(\sigma) + o(d_*^{-2})$. Finally, we have $\mathrm{E}|\hat{\eta}|\mathbf{1}_{\mathcal{E}^c} = o(d_*^{-2})$.

5.9. *Regarding Section* 4.1. It suffices to show that (3) and (4) hold. Let $i'$, $j'$, $k'$ be a permutation of $i$, $j$, $k$ and w.l.o.g. let $d_{i'} = d_{i'} \wedge d_{j'} \wedge d_{k'} = d_i \wedge d_j \wedge d_k$. Then by Lemmas 5.1 and 5.2, it can be shown that

$$\mathrm{E}\left|\frac{1}{\sqrt{d_j d_k}}[y'PV_{i'}PV_{j'}PV_{k'}Py - \mathrm{E}(y'PV_{i'}PV_{j'}PV_{k'}Py)]\right|^g \le c.$$

Here $c$ represents a constant whose value may be different at different places. Similarly, it can be shown that $|\operatorname{tr}(PV_{i'}PV_{j'}PV_{k'})| \leq c(d_i d_j d_k / d_i \vee d_j \vee d_k)$,

$$\mathrm{E}\left(\frac{d_*^2}{d_i d_j d_k d_l} \sup_{\tilde{\sigma} \in S_\delta} |y'\tilde{P}V_{i'}\tilde{P}V_{j'}\tilde{P}V_{k'}\tilde{P}V_{l'}Py|\right)^g \leq c$$

and $(d_*^2/d_i d_j d_k d_l)|\operatorname{tr}(\tilde{P}V_{i'}\tilde{P}V_{j'}\tilde{P}V_{k'}\tilde{P}V_{l'})| \leq c$.

As for (4), first note that $P(\sigma) = -(1/2)(\operatorname{tr}(HG_i HG_j)/d_i d_j)_{0 \leq i,j \leq q}$. It can be shown that

$$|\operatorname{tr}(\tilde{H}G_i \tilde{H}G_j) - \operatorname{tr}(HG_i HG_j)|$$

$$\leq \sum_k |\tilde{\sigma}_k - \sigma_k||\operatorname{tr}(HG_k \tilde{H}G_i HG_j)|$$

$$+ \sum_l |\tilde{\sigma}_l - \sigma_l||\operatorname{tr}(HG_l \tilde{H}G_j HG_i)|$$

$$+ \sum_{k,l} |\tilde{\sigma}_k - \sigma_k||\tilde{\sigma}_l - \sigma_l||\operatorname{tr}(HG_k \tilde{H}G_i HG_l \tilde{H}G_j)|,$$

$|\operatorname{tr}(HG_k \tilde{H}G_i HG_j)| \leq 2\sigma_k^{-1} d_i d_j$ and $|\operatorname{tr}(HG_k \tilde{H}G_i HG_l \tilde{H}G_j)| \leq 4(\sigma_k \sigma_l)^{-1} \times d_i d_j$. Thus $\sup_{\tilde{\sigma} \in S_\delta} \|P(\tilde{\sigma}) - P(\sigma)\| \to 0$ as $\delta \to 0$. Note that $U_1(\sigma) = -P(\sigma)$. Similarly, one can show that $\sup_{\tilde{\sigma} \in S_\delta} \|U_2(\tilde{\sigma}) - U_2(\sigma)\| \to 0$ as $\delta \to 0$. Finally, it can be shown that

$$u_{ijklm} = \frac{1}{2d_j d_k d_l d_m}\bigg\{[\operatorname{tr}(PV_i PV_j PV_l) + \operatorname{tr}(PV_i PV_l PV_j)]\operatorname{tr}(PV_k PV_m)$$

$$+ \sum_{a,b,c} \operatorname{tr}(PV_a PV_b PV_c PV_k PV_m)$$

$$+ \sum_{a,b,c} \operatorname{tr}(PV_a PV_b PV_c PV_m PV_k)\bigg\},$$

where the summation is over all $a$, $b$, $c$ which is a permutation of $i$, $j$, $l$. It follows that $\sup_{\tilde{\sigma} \in S_\delta} \|U_3(\tilde{\sigma}) - U_3(\sigma)\| \to 0$ as $\delta \to 0$.

5.10. *Regarding Section* 4.2. First note that the formulas derived in Sections 2.1 and 5.2 for $l_{\mathrm{R}}(\sigma)$ and its derivatives hold for the general linear mixed model (1.1), which includes the longitudinal model. Next, note that the matrix $P$ of (2.6) can be expressed as $P = \Sigma^{-1} + \Delta$, where $\|\Delta\|_2$ is bounded. These results and Lemma 5.1 are used to verify the moment conditions involved.

As for the conditions regarding $\eta$ and its derivatives, we have the following expressions: $\eta(\sigma) = g_1(\sigma) + g_2(\sigma) + g_3(\sigma)$, where

$$g_1(\sigma) = \sum_{i=1}^{t} m_i'(G_i - G_i Z_i' \Sigma_i^{-1} Z_i G_i) m_i,$$

$$g_2(\sigma) = \left( l - \sum_{i=1}^{t} X_i' \Sigma_i^{-1} Z_i G_i m_i \right)' \left( \sum_{i=1}^{t} X_i' \Sigma_i^{-1} X_i \right)^{-1} \left( l - \sum_{i=1}^{t} X_i' \Sigma_i^{-1} Z_i G_i m_i \right),$$

$$g_3(\sigma) = \sum_{i=1}^{t} \mathrm{tr} \left\{ \left[ \frac{\partial}{\partial \sigma} (\Sigma_i^{-1} Z_i G_i m_i) \right]' \Sigma_i \left[ \frac{\partial}{\partial \sigma} (\Sigma_i^{-1} Z_i G_i m_i) \right] \right\}.$$

With these expressions one can verify the conditions regarding $\eta$ $b$ and $\|B\|$.

Finally, for condition (iv) of Theorem 4.1 we have, for example, $\tilde{Q} - Q = -(1/2t)[\mathrm{tr}(P\Sigma_j P\Sigma_k)]_{0 \leq j,k \leq q}$, where $\Sigma_j = \partial\Sigma/\partial\sigma_j$. Note that $\tilde{P} - P = P(\Sigma - \tilde{\Sigma})\tilde{P}$, where $\tilde{P}$ is $P$ with $\sigma$ replaced by $\tilde{\sigma}$ and so on.

K. Das
Department of Statistics
Calcutta University
35 Ballygunge Circular Road
Calcutta 700 019
India
e-mail: kalyan_stat@Yahoo.co.in

J. Jiang
Department of Statistics
University of California
One Shields Avenue
Davis, California 95616
USA
e-mail: jiang@wald.ucdavis.edu

J. N. K. Rao
School of Mathematics
  and Statistics
Carleton University
Ottawa, Ontario
Canada K1S 5B6
e-mail: jrao@math.carleton.ca